\documentclass[12pt]{article}
\usepackage{latexsym,amssymb,amsmath,enumerate,float,geometry,cite,algorithm2e}
\newtheorem{theorem}{Theorem}

\newtheorem{claim}{Claim}

\begin{document}

\title{Dynamic Monopolies for Degree Proportional Thresholds\\ in Connected Graphs of Girth at least Five and Trees}
\author{Michael Gentner and Dieter Rautenbach\\
Institute of Optimization and Operations Research\\Ulm University, Ulm, Germany\\
$\{$michael.gentner,dieter.rautenbach$\}$@uni-ulm.de}
\date{}
\maketitle

\begin{abstract}
Let $G$ be a graph, and let $\rho\in [0,1]$.
For a set $D$ of vertices of $G$, let the set $H_{\rho}(D)$ arise by starting with the set $D$, 
and iteratively adding further vertices $u$ to the current set if they have at least $\lceil \rho d_G(u)\rceil$ neighbors in it. 
If $H_{\rho}(D)$ contains all vertices of $G$, then $D$ is known as an irreversible dynamic monopoly or a perfect target set
associated with the threshold function $u\mapsto \lceil \rho d_G(u)\rceil$.
Let $h_{\rho}(G)$ be the minimum cardinality of such an irreversible dynamic monopoly.

For a connected graph $G$ of maximum degree at least $\frac{1}{\rho}$, 
Chang (Triggering cascades on undirected connected graphs, Information Processing Letters 111 (2011) 973-978) showed
$h_{\rho}(G)\leq 5.83\rho n(G)$,
which was improved by 
Chang and Lyuu (Triggering cascades on strongly connected directed graphs, Theoretical Computer Science 593 (2015) 62-69) to 
$h_{\rho}(G)\leq 4.92\rho n(G)$.
We show that for every $\epsilon>0$, there is some $\rho(\epsilon)>0$ such that 
$h_{\rho}(G) \leq(2+\epsilon)\rho n(G)$
for every $\rho$ in $(0,\rho(\epsilon))$, 
and every connected graph $G$ 
that has maximum degree at least $\frac{1}{\rho}$ and girth at least $5$.
Furthermore, we show that 
$h_{\rho}(T) \leq \rho n(T)$
for every $\rho$ in $(0,1]$, 
and every tree $T$ that has order at least $\frac{1}{\rho}$.
\end{abstract}

{\small
\begin{tabular}{lp{12.5cm}}
\textbf{Keywords:} & Irreversible dynamic monopoly; perfect target set\\
\textbf{MSC2010:} & 05C69
\end{tabular}
}

\section{Introduction}

We consider finite, simple, and undirected graphs, and use standard terminology and notation.

Let $G$ be a graph with vertex set $V(G)$.
Let $\phi:V(G)\to \mathbb{N}_0$ be a threshold function such that $\phi(u)$ is at most the degree $d_G(u)$ of $u$ in $G$ for every vertex $u$ of $G$.
For a set $D$ of vertices of $G$, let $H_{(G,\phi)}(D)$ be the smallest set $\bar{D}$ of vertices of $G$ such that $D\subseteq \bar{D}$,
and every vertex $u$ in $V(G)\setminus \bar{D}$ has less than $\phi(u)$ neighbors in $\bar{D}$.
Note that the set $H_{(G,\phi)}(D)$ can be constructed by starting with the set $D$, 
and iteratively adding further vertices $u$ to the current set if they have at least $\phi(u)$ neighbors in it. 
Such iterative expansion processes have been considered in a variety of contexts \cite{abtz,cdprs,cprs,cr,c,dr,ksz,p,z}.
If $H_{(G,\phi)}(D)=V(G)$, then $D$ is a {\it $\phi$-dynamic monopoly} of $G$.
Let $h_{\phi}(G)$ be the minimum cardinality of a $\phi$-dynamic monopoly of $G$.

By a simple probabilistic argument, 
very similar to the one used by Alon and Spencer \cite{as} to prove the Caro-Wei bound on the independence number of a graph \cite{c,w}, 
Ackerman, Ben-Zwi, Wolfovitz \cite{abw} showed
\begin{eqnarray}\label{e1}
h_{\phi}(G) & \leq & \sum_{u\in V(G)}\frac{\phi(u)}{d_G(u)+1}.
\end{eqnarray}
Essentially the same bound was obtained independently by Reichman \cite{r}.
For an application of this argument to independence in hypergraphs see \cite{bghr}.

It is easy to see \cite{cr} that (minimum) dynamic monopolies and (maximum) generalized degenerate induced subgraphs are dual notions,
that is, (\ref{e1}) can be considered the dual counterpart of bounds as in \cite{aks}.

A very natural choice for the threshold function is to assign values that are proportional to the vertex degrees.
Specifically, for some real parameter $\rho$ in $[0,1]$, let 
$$\phi_{\rho}:V(G)\to \mathbb{N}_0:u\mapsto \lceil \rho d_G(u)\rceil.$$
If $D$ is a random set of vertices of $G$ that contains each vertex independently at random with probability $\rho$,
then, for every vertex $u$ of $G$, the expected number of neighbors of $u$ that belong to $D$ is $\rho d_G(u)$.
This suggests that $h_{\phi_{\rho}}(G)$ might be only slightly bigger than $\rho n(G)$, where $n(G)$ is the order of $G$.
Without any restriction on $\rho$ or $G$ though, this intuition is misleading. 
In fact, if $\rho$ is positive but much smaller than $\frac{1}{n(G)}$, 
then $h_{\phi_{\rho}}(G)$ is at least $1$, while $\rho n(G)$ can be arbitrarily small.  
As observed by Chang \cite{c2}, it is reasonable to consider only values of $\rho$ that are at least $\frac{1}{\Delta(G)}$,
where $\Delta(G)$ is the maximum degree of $G$,
because $\phi_{\frac{1}{\Delta(G)}}=\phi_{\rho}$ for every $\rho$ in $\left(0,\frac{1}{\Delta(G)}\right]$.
For a connected graph $G$ and $\rho\in \left[\frac{1}{\Delta(G)},1\right]$,
Chang \cite{c2} proved
\begin{eqnarray}\label{e2}
h_{\phi_{\rho}}(G) & \leq & (2\sqrt{2}+3)\rho n(G) \approx 5.83 \rho n(G),
\end{eqnarray}
which was improved by Chang and Lyuu \cite{cl} to 
\begin{eqnarray}\label{e3}
h_{\phi_{\rho}}(G) & \leq & 4.92\rho n(G).
\end{eqnarray}
Note that the bound in (\ref{e1}) might evaluate to $\Omega(n(G))$, 
because, for instance, vertices of degree $1$ contribute $\frac{1}{2}$ rather than $O(\rho)$ to the right hand side of (\ref{e1}).
In fact, especially for small values of $\rho$, and graphs with many vertices of small degrees,
the bound (\ref{e3}) can be much better than the bound (\ref{e1}).

The proof strategies for (\ref{e2}) and (\ref{e3}) are quite different.
The bound (\ref{e2}) is proved by a suitable adaptation of the argument of Ackerman, Ben-Zwi, Wolfovitz \cite{abw}.
Vertices of small degree, that is, at most $\frac{1}{\rho}$, are treated differently from those of large degree,
that is, more than $\frac{1}{\rho}$.
A small set $X_0$ of vertices of large degree ensures 
that the remaining vertices of large degree have few neighbors of small degree outside of $H_{(G,\phi_{\rho})}(X_0)$.
This allows to apply the argument of Ackerman et al. to the vertices of large degree outside of $X_0$.
The bound (\ref{e3}) is proved by a random procedure that considers a sequence $X_1,X_2,\ldots$ of random sets of vertices 
each containing every individual vertex independently at random with probability $3.51\rho$.
Starting with the empty set, a $\phi_{\rho}$-dynamic monopoly is constructed by iteratively adding
the vertices in $X_i\setminus H_{(G,\phi_{\rho})}(X_1\cup \ldots \cup X_{i-1})$ to the current set.
Chernoff's inequality is used to ensure that $H_{(G,\phi_{\rho})}(X_1\cup \ldots \cup X_i)$ grows sufficiently fast.
The proof of (\ref{e3}) has some resemblance to iterative random procedures 
that are used to show lower bounds on the independence number \cite{ghrs,lw}.

It is a natural to ask for the best-possible constant in bounds of the form (\ref{e2}) and (\ref{e3}).
We contribute to this question by showing the following results.

\begin{theorem}\label{theorem1}
For every positive real $\epsilon$, there is some positive real $\rho(\epsilon)$ such that 
\begin{eqnarray}\label{e4}
h_{\phi_{\rho}}(G) & \leq & (2+\epsilon)\rho n(G)
\end{eqnarray}
for every $\rho$ in $(0,\rho(\epsilon))$, 
and every connected graph $G$ 
that has maximum degree at least $\frac{1}{\rho}$ and girth at least $5$.
\end{theorem}
The proof of Theorem \ref{theorem1} is based on a combination of the techniques from \cite{c2,cl}.
Note that (\ref{e4}) requires a sufficiently small value of $\rho$,
but that bounds like (\ref{e2}), (\ref{e3}), and (\ref{e4}) are especially interesting for small values of $\rho$.
It is possible to generalize (\ref{e4}) to strongly connected directed graphs similarly as in \cite{cl}.

\begin{theorem}\label{theorem2}
If $\rho$ is in $(0,1]$, and $T$ is a tree of order at least $\frac{1}{\rho}$, then 
\begin{eqnarray}\label{e4}
h_{\phi_{\rho}}(T) & \leq & \rho n(T).
\end{eqnarray}
\end{theorem}
Note that $h_{\phi_{\rho}}(T)$ can be computed in linear time \cite{cdprs} for a given tree.

For many more references to and discussion of related work see \cite{c2,cl}.

\section{Proofs of Theorem \ref{theorem1} and Theorem \ref{theorem2}}

We recall two tools from probability theory.
\begin{itemize}
\item {\bf Markov's inequality} (cf. Theorem 3.2 in \cite{mr})\\
If $X$ is a non-negative random variable, and $t>0$, then $\mathbb{P}\big[X\geq t\big]\leq \frac{1}{t}\mathbb{E}\big[X\big]$.
\item {\bf Chernoff's inequality} (cf. Theorem 4.2 in \cite{mr})\\
If $Z_1,\ldots,Z_n$ are independent random variables, and $p_1,\ldots,p_n$ in $(0,1)$ are such that 
$\mathbb{P}[Z_i=1]=p_i$ and $\mathbb{P}[Z_i=0]=1-p_i$ for $i\in [n]$, then
$$\mathbb{P}\Big[Z_1+\cdots+Z_n<(1-\delta)\mathbb{E}\big[Z_1+\cdots+Z_n\big]\Big]
<e^{-\frac{\delta^2}{2}\mathbb{E}\big[Z_1+\cdots+Z_n\big]}$$
for every $0<\delta\leq 1$.
\end{itemize}
We proceed to the first proof.

\medskip

\noindent {\it Proof of Theorem \ref{theorem1}:}
Let $\epsilon>0$.
Let $\delta>0$ be small enough such that
\begin{eqnarray}\label{edelta}
\delta & \leq & \min\left\{ e^{-\frac{1}{4}},\frac{1}{2}\right\}\mbox{, and}\label{edelta1}\\
(1+\delta)^2+\frac{1+\delta}{(1-\delta)^2} & \leq & 2+\epsilon.\label{edelta2}
\end{eqnarray}
Let $\rho(\epsilon)>0$ be small enough such that 
\begin{eqnarray}
\rho(\epsilon) & \leq & \frac{\delta}{1+\delta}\left(1-e^{-\frac{\delta^2}{2(1-\delta)}}\right)\frac{1}{8\ln\left(\frac{1}{\delta}\right)}.\label{erho}
\end{eqnarray}
Let $\rho\in (0,\rho(\epsilon))$. Note that (\ref{erho}) implies $\rho\leq \delta$.
Let $G$ be a connected graph of maximum degree at least $\frac{1}{\rho}$ and girth at least $5$.

Let
\begin{eqnarray*}
V_1&=&\left\{ u\in V(G):d_G(u)<\frac{1}{\rho}\right\}\mbox{, and }\\
V_2&=&V(G)\setminus V_1.
\end{eqnarray*}
Note that $\phi_{\rho}(u)=1$ for $u\in V_1$, and $V_2\not=\emptyset$.

Throughout this proof, we write `$H(X)$' instead of `$H_{(G,\phi_{\rho})}(X)$'.

\begin{claim}\label{claim1}
There is a set $X_0\subseteq V_2$ with the following properties. 
\begin{enumerate}[(i)]
\item $|N_G(u)\cap (V_1\setminus H(X_0))|\leq \frac{1}{1+\delta}d_G(u)$ for every $u\in V_2\setminus X_0$.
\item $|X_0|\leq (1+\delta)\rho n(G)$.
\end{enumerate}
\end{claim}
{\it Proof:} Let $(x_1,\ldots,x_k)$ be a maximal sequence of distinct vertices from $V_2$ such that 
$$\Big|N_G(x_i)\cap \Big(V_1\setminus H(\{ x_j:1\leq j\leq i-1\})\Big)\Big|>\frac{1}{1+\delta}d_G(x_i)$$ 
for $1\leq i\leq k$.

Let $X_0=\{ x_1,\ldots,x_k\}$.
The maximality of the sequence $(x_1,\ldots,x_k)$ implies (i).

Since $\phi_{\rho}(u)=1$ for $u\in V_1$, and $d_G(x_i)\geq \frac{1}{\rho}$ for $1\leq i\leq k$, 
we have 
\begin{eqnarray*}
|H(\{ x_j:1\leq j\leq i\})|-|H(\{ x_j:1\leq j\leq i-1\})| & \geq &
\left| N_G(x_i)\cap \Big(V_1\setminus H(\{ x_j:1\leq j\leq i-1\})\Big)\right|\\
& > & \frac{1}{1+\delta}d_G(x_i)\\
& \geq & \frac{1}{(1+\delta)\rho},
\end{eqnarray*}
which implies 
$\frac{k}{(1+\delta)\rho}<|H(\{ x_1,\ldots,x_k\})|\leq n(G)$.
It follows that $|X_0|=k\leq (1+\delta)\rho n(G)$, which completes the proof of the claim.
$\Box$

\medskip

\noindent For $i\in \mathbb{N}$, let $X_i$ be a random subset of $V(G)\setminus H(X_0)$ 
that contains each vertex of $V(G)\setminus H(X_0)$ independently at random with probability 
\begin{eqnarray}\label{e5}
p^{(1)}=\frac{\rho}{(1-\delta)}.
\end{eqnarray}
Note that $\rho\leq \delta$ implies $p^{(1)}\leq \frac{\delta}{(1-\delta)}\stackrel{(\ref{edelta})}{\leq} 1$.
Furthermore, note that $X_0$ is chosen deterministically,
and that the random sets $X_1,X_2,\ldots$ are chosen independently of each other.

Let $Y_0=X_0$, and, for $i\in \mathbb{N}$, let
\begin{eqnarray*}
Y_i & = & X_i\setminus H(Y_0\cup \ldots \cup Y_{i-1})\mbox{, and}\\
Y_{\leq i} & = & Y_0\cup Y_1\cup \ldots \cup Y_i.
\end{eqnarray*}
By construction, 
$X_0\cup X_1\cup \ldots \cup X_i\subseteq H(Y_{\leq i-1})\cup Y_i\subseteq H(Y_{\leq i})$,
which implies 
$$H(X_0\cup X_1\cup \ldots \cup X_i)=H(Y_{\leq i}).$$
\begin{claim}\label{claim2}
$\mathbb{P}\big[u\not\in H(X_0\cup X_i)\big]\leq \delta$
for $i\in \mathbb{N}$ and $u\in V_2\setminus H(X_0)$.
\end{claim}
{\it Proof:} Let 
\begin{eqnarray*}
N_1 & = & N_G(u)\cap (V_1\setminus H(X_0)),\\
N_2 & = & N_G(u)\setminus (V_1\cup H(X_0))\mbox{, and }\\ 
N_3 & = & N_G(u)\cap H(X_0).
\end{eqnarray*}
By Claim \ref{claim1}, $|N_1|\leq \frac{1}{1+\delta}d_G(u)$, 
which implies 
\begin{eqnarray}\label{en2n3}
|N_2|+|N_3| & = & d_G(u)-|N_1|\geq \frac{\delta}{1+\delta}d_G(u)\stackrel{u\in V_2}{\geq}\frac{\delta}{(1+\delta)\rho}.
\end{eqnarray}
Let $G_u$ be the subgraph of $G$ that arises from the induced subgraph 
$$G\left[N_G(u)\cup \bigcup_{v\in N_2}\left(N_G(v)\setminus \{ u\}\right)\right]$$
by removing, for every vertex $v$ in $N_2$ and every vertex $w$ in $N_G(v)\setminus \{ u\}$,
every edge incident with $w$ except for the edge $vw$.
Since $G$ has girth at least $5$, 
all vertices in $N_1\cup N_3$ are isolated in $G_u$,
and each vertex $v$ in $N_2$ belongs to a component of $G_u$ that is a star of order $d_G(v)$ with center $v$
whose set of endvertices is $N_G(v)\setminus \{ u\}$.

Let 
$$H_u=H_{(G_u,\phi_{\rho})}\Big((H(X_0)\cup X_i)\cap V(G_u)\Big).$$
Since $G_u$ is a subgraph of $G$, we have $H_u\subseteq H(H(X_0)\cup X_i)=H(X_0\cup X_i)$.

The structure of $G_u$ and the choice of $X_i$ implies that 
\begin{itemize}
\item every vertex in $N_1$ belongs to $H_u$ with probability $p^{(1)}$, and 
\item every vertex in $N_3$ belongs to $H_u$ with probability $1$.
\end{itemize}
Let $v\in N_2$.
Note that $v\not\in H_u$ implies 
$v\not\in X_i$ as well as $|(N_G(v)\setminus \{ u\})\cap X_i|<\rho d_G(v)$,
that is,
$|(N_G[v]\setminus \{ u\})\cap X_i|<\rho d_G(v)$,
where $N_G[v]=\{ v\}\cup N_G(v)$.
Since $|N_G[v]\setminus \{ u\}|=d_G(v)$, 
the definitions of $p^{(1)}$ and $X_i$ imply
\begin{eqnarray*}
\mathbb{P}\Big[v\not\in H_u\Big] & \leq & \mathbb{P}\Big[|(N_G[v]\setminus \{ u\})\cap X_i|<\rho d_G(v)\Big]\\ 
& \stackrel{(\ref{e5})}{=} & \mathbb{P}\left[|(N_G[v]\setminus \{ u\})\cap X_i|<(1-\delta)p^{(1)}d_G(v)\right]\\
& = & \mathbb{P}\Big[|(N_G[v]\setminus \{ u\})\cap X_i|<(1-\delta)\mathbb{E}\big[|(N_G[v]\setminus \{ u\})\cap X_i|\big]\Big]\\
& < & e^{-\delta^2\mathbb{E}\big[|(N_G[v]\setminus \{ u\})\cap X_i|\big]/2},
\end{eqnarray*}
where the last inequality is a consequence of Chernoff's inequality.

Since 
$$\mathbb{E}\big[|(N_G[v]\setminus \{ u\})\cap X_i|\big]
=p^{(1)}d_G(v)
\stackrel{(\ref{e5})}{=}
\frac{\rho d_G(v)}{(1-\delta)}
\stackrel{v\in V_2}{\geq}\frac{1}{(1-\delta)},$$
we obtain
$$\mathbb{P}\big[v\not\in H_u\big] < e^{-\frac{\delta^2}{2(1-\delta)}}.$$
This implies that
\begin{itemize}
\item every vertex in $N_2\cup N_3$ belongs to $H_u$ with probability more than $p^{(2)}$, where
\end{itemize}
\begin{eqnarray}\label{ep2}
p^{(2)}=1-e^{-\frac{\delta^2}{2(1-\delta)}}.
\end{eqnarray}
By the structure of $G_u$, the events $\big[ v\in H_u\big]$ are all independent for $v\in N_G(u)$.
Let $X_u$ be a random subset of $N_2\cup N_3$ that contains each vertex of $N_2\cup N_3$ independently at random with probability $p^{(2)}$.

Note that 
\begin{eqnarray}\label{eebig}
\mathbb{E}\Big[|(N_2\cup N_3)\cap X_u|\Big]=p^{(2)}(|N_2|+|N_3|)
\stackrel{(\ref{en2n3})}{\geq} \frac{\delta p^{(2)}}{(1+\delta)\rho}
\stackrel{(\ref{erho})}{\geq} 8\ln\left(\frac{1}{\delta}\right),
\end{eqnarray}
and
\begin{eqnarray}\label{eepssm}
\frac{\rho d_G(u)}{p^{(2)}(|N_2|+|N_3|)}
\stackrel{(\ref{en2n3})}{\leq}\frac{\rho (1+\delta)}{\delta p^{(2)}}
\stackrel{(\ref{erho})}{\leq} \frac{1}{8\ln\left(\frac{1}{\delta}\right)}
\stackrel{(\ref{edelta1})}{\leq} \frac{1}{2}.
\end{eqnarray}
Since $u\not\in H(X_0\cup X_i)$ implies $|(N_2\cup N_3)\cap H_u|<\rho d_G(u)$,
the definition of $X_u$ and Chernoff's inequality imply
\begin{eqnarray*}
\mathbb{P}\Big[u\not\in H(X_0\cup X_i)\Big]
& \leq & \mathbb{P}\Big[|(N_2\cup N_3)\cap H_u|<\rho d_G(u)\Big]\\
& \leq & \mathbb{P}\Big[|(N_2\cup N_3)\cap X_u|<\rho d_G(u)\Big]\\
& = & \mathbb{P}\left[|(N_2\cup N_3)\cap X_u|<\frac{\rho d_G(u)}{p^{(2)}(|N_2|+|N_3|)}p^{(2)}(|N_2|+|N_3|)\right]\\
& \stackrel{(\ref{eebig}),(\ref{eepssm})}{\leq} & \mathbb{P}\left[|(N_2\cup N_3)\cap X_u|<\frac{1}{2}\mathbb{E}\Big[|(N_2\cup N_3)\cap X_u|\Big]\right]\\
& < & e^{-\frac{1}{8}\mathbb{E}\big[|(N_2\cup N_3)\cap X_u|\big]}\\
& \stackrel{(\ref{eebig})}{\leq} & \delta.
\end{eqnarray*}
$\Box$

\begin{claim}\label{claim3}
$\mathbb{P}\big[u\not\in H(Y_{\leq i})\mid u\not\in H(Y_{\leq i-1})\big]\leq\delta$
for $i\in \mathbb{N}$ and $u\in V(G)\setminus H(X_0)$.
\end{claim}
{\it Proof:} 
Recall that the set $X_0$ is chosen deterministically,
that is, the only source of randomness are the sets $X_1,X_2,\ldots$, 
which are chosen independently of each other.
This implies that for every two not necessarily distinct vertices $u'$ and $u''$ of $V(G)$,
the two events $\big[u'\not\in H(X_0\cup X_i)\big]$ and $\big[u''\not\in H(Y_{\leq i-1})\big]$ are independent.

First, let $u\in V_2$.
Since $H(X_0\cup X_i)\subseteq H(Y_{\leq i})$, 
Claim \ref{claim2} and the independence observed above imply
\begin{eqnarray*}
\mathbb{P}\big[u\not\in H(Y_{\leq i})\mid u\not\in H(Y_{\leq i-1})\big]
&\leq & \mathbb{P}\big[u\not\in H(X_0\cup X_i)\mid u\not\in H(Y_{\leq i-1})\big]\\
&= & \mathbb{P}\big[u\not\in H(X_0\cup X_i)\big]\\
& \leq & \delta.
\end{eqnarray*}
Next, let $u\in V_1$.
Let $u_0\ldots u_{\ell}$ be a shortest path in $G$ between $u=u_0$ and some vertex $u_{\ell}$ in $V_2$.
Since $\phi_{\rho}(u_0)=\ldots=\phi_{\rho}(u_{\ell-1})=1$,
we obtain that $u\not\in H(Y_{\leq i})$ implies $u_{\ell}\not\in H(Y_{\leq i})$,
which implies $u_{\ell}\not\in H(X_0\cup X_i)$.
Now, Claim \ref{claim2} and the independence observed above imply
\begin{eqnarray*}
\mathbb{P}\big[u\not\in H(Y_{\leq i})\mid u\not\in H(Y_{\leq i-1})\big]
& \leq & \mathbb{P}\big[u_{\ell}\not\in H(Y_{\leq i})\mid u\not\in H(Y_{\leq i-1})\big]\\
& \leq & \mathbb{P}\big[u_{\ell}\not\in H(X_0\cup X_i)\mid u\not\in H(Y_{\leq i-1})\big]\\
& = & \mathbb{P}\big[u_{\ell}\not\in H(X_0\cup X_i)\big]\\
& \leq & \delta.
\end{eqnarray*}
$\Box$

\medskip

\noindent Note that $u\not\in H(Y_{\leq i})$ implies $u\not\in H(Y_{\leq i-1})$,
and that $\mathbb{P}\big[u\not\in H(Y_{\leq i})\big]=0$ for $u\in H(X_0)$.
By Claim \ref{claim3}, for $i\in \mathbb{N}$, we obtain, by linearity of expectation,
\begin{eqnarray}
\mathbb{E}\big[\big|V(G)\setminus H(Y_{\leq i})\big|\big]
&=& \sum_{u\in V(G)\setminus H(X_0)}\mathbb{P}\big[u\not\in H(Y_{\leq i})\big]\nonumber\\
&=& \sum_{u\in V(G)\setminus H(X_0)}\mathbb{P}\big[u\not\in H(Y_{\leq i})\mid u\not\in H(Y_{\leq i-1})|\big]\mathbb{P}\big[u\not\in H(Y_{\leq i-1})\big]\nonumber\\
&\leq & \delta\sum_{u\in V(G)\setminus H(X_0)}\mathbb{P}\big[u\not\in H(Y_{\leq i-1})\big]\nonumber\\
& = & \delta\mathbb{E}\big[\big|V(G)\setminus H(Y_{\leq i-1})\big|\big],\label{ey}
\end{eqnarray}
which is equivalent to 
\begin{eqnarray}\label{ee}
\mathbb{E}\big[\big|V(G)\setminus H(Y_{\leq i-1})\big|\big] \leq 
\frac{1}{1-\delta}\left(\mathbb{E}\big[\big|H(Y_{\leq i})\big|\big]-\mathbb{E}\big[\big|H(Y_{\leq i-1})\big|\big]\right).
\end{eqnarray}
Iteratively applying (\ref{ey}), we obtain
\begin{eqnarray*}
\mathbb{E}\big[\big|V(G)\setminus H(Y_{\leq i})\big|\big]
& \leq & \delta^i n(G).
\end{eqnarray*}
Furthermore, by the choice of $X_i$ and the definition of $Y_i$, we have 
$\mathbb{P}\big[u\in Y_i\mid u\not\in H(Y_{\leq i-1})\big]=p^{(1)}$
and
$\mathbb{P}\big[u\in Y_i\mid u\in H(Y_{\leq i-1})\big]=0$, which implies
\begin{eqnarray*}
\mathbb{E}\big[|Y_i|\big]
&=& \sum_{u\in V(G)\setminus H(X_0)}\mathbb{P}\big[u\in Y_i\big]\\
&=& \sum_{u\in V(G)\setminus H(X_0)}\mathbb{P}\big[u\in Y_i\mid u\not\in H(Y_{\leq i-1})\big]\mathbb{P}\big[u\not\in H(Y_{\leq i-1})\big]\\
&=& p^{(1)}\sum_{u\in V(G)\setminus H(X_0)}\mathbb{P}\big[u\not\in H(Y_{\leq i-1})|\big]\\
&=& p^{(1)}\mathbb{E}\big[\big|V(G)\setminus H(Y_{\leq i-1})\big|\big]\\
& \stackrel{(\ref{ee})}{\leq} & \frac{p^{(1)}}{1-\delta}\left(\mathbb{E}\big[\big|H(Y_{\leq i})\big|\big]-\mathbb{E}\big[\big|H(Y_{\leq i-1})\big|\big]\right).
\end{eqnarray*}
Therefore, by Claim \ref{claim1},
\begin{eqnarray*}
\mathbb{E}\big[|Y_{\leq i}|\big] & = & \mathbb{E}\big[|Y_0\cup Y_1\cup \ldots \cup Y_i|\big]\\ 
& \leq & |X_0|+\mathbb{E}\big[|Y_1\cup\ldots \cup Y_i|\big]\\
& \leq & |X_0|+\frac{p^{(1)}}{1-\delta}\sum_{j=1}^{i}
\left(\mathbb{E}\big[\big|H(Y_{\leq j})\big|\big]-\mathbb{E}\big[\big|H(Y_{\leq j-1})\big|\big]\right)\\
& \leq & |X_0|+\frac{p^{(1)}}{1-\delta}n(G)\\
& \stackrel{(\ref{e5})}{=} & |X_0|+\frac{\rho}{(1-\delta)^2}n(G)\\
& \leq & \rho(1+\delta)n(G)+\frac{\rho}{(1-\delta)^2}n(G)\\
& = & \left((1+\delta)+\frac{1}{(1-\delta)^2}\right)\rho n(G).
\end{eqnarray*}
Let $k\in \mathbb{N}$ be large enough such that $\delta^k n(G)+\frac{1}{1+\delta}<1$.

By Markov's inequality,
\begin{eqnarray*}
&& \mathbb{P}\left[\left(\big|V(G)\setminus H(Y_{\leq k})\big|<1\right)\wedge \left(|Y_{\leq k}|<(1+\delta)\mathbb{E}\big[|Y_{\leq k}|\big]\right)\right]\\
& \geq & 1
-\mathbb{P}\left[\big|V(G)\setminus H(Y_{\leq k})\big|\geq 1\right]
-\mathbb{P}\left[|Y_{\leq k}|\geq (1+\delta)\mathbb{E}\big[|Y_{\leq k}|\big]\right]\\
& \geq & 1-\delta^k n(G)-\frac{1}{1+\delta}\\
& > & 0.
\end{eqnarray*}
By the first moment method \cite{as}, 
this implies the existence of a set $D$ of vertices such that
\begin{eqnarray*}
|D| & \leq & (1+\delta)\mathbb{E}\big[|Y_{\leq k}|\big]\\
& \leq & (1+\delta)\left((1+\delta)+\frac{1}{(1-\delta)^2}\right)\rho n(G)\\
& \stackrel{(\ref{edelta2})}{\leq} & (2+\epsilon)\rho n(G),
\end{eqnarray*}
and 
$H(D)=V(G)$,
which completes the proof of Theorem \ref{theorem1}.
$\Box$

\medskip

\noindent It is conceivable that (\ref{e1}) remains true 
without the girth condition
as well as for arbitrary values of $\rho$ in $\left[\frac{1}{\Delta(G)},1\right]$ and sufficiently large girth.

We proceed to the second proof.

\medskip

\noindent {\it Proof of Theorem \ref{theorem2}:}
Let $\rho$ be in $(0,1]$.
Suppose that the tree $T$ is a counterexample of minimum order,
that is, $n(T)\geq \frac{1}{\rho}$ and $h_{\phi_{\rho}}(T) > \rho n(T)$.
Let $V_2=\left\{ u\in V(T):d_T(u)\geq \frac{1}{\rho}\right\}$.
Note that $V_2$ is a $\phi_{\rho}$-dynamic monopoly of $T$.
If $|V_2|\leq 1$, then $h_{\phi_{\rho}}(T)\leq |V_2|\leq \rho n(T)$, which is a contradiction.
Hence, $|V_2|\geq 2$.

Let $u$ in $V_2$ be chosen such that a largest component $K$ of $T-u$ that contains a vertex from $V_2$ has largest possible order.
Suppose that $T-u$ has a second component $K'$ distinct from $K$ that contains a vertex from $V_2$.
Let $u'$ be in $V(K')\cap V_2$.
Since $T-u'$ has a component that contains all vertices in $\{ u\}\cup V(K)$,
we obtain a contradiction to the choice of $u$.
Hence, $K$ is the only component of $T-u$ that contains a vertex from $V_2$.
Let $R=T-V(K)$.
Since $K$ and $R$ both contain a vertex from $V_2$, we have $n(K)\geq \frac{1}{\rho}$ and $n(R)\geq \frac{1}{\rho}$.
By the choice of $T$, the tree $K$ has a $\phi_{\rho}$-dynamic monopoly $D$ with $|D|\leq \rho n(K)$.
Let $v$ be the neighbor of $u$ in $K$.
Since $\rho d_T(v)\leq 1+\rho(d_T(v)-1)=1+\rho d_K(v)$, the set $\{ u\}\cup D$ is a $\phi_{\rho}$-dynamic monopoly of $T$ 
of order at most $|D|+1\leq \rho n(K)+1\leq \rho n(K)+\rho n(R)=\rho n(T)$.
This contradiction completes the proof. $\Box$

\medskip

\noindent If $\rho\in (0,1]$ is such that $\frac{1}{\rho}$ is an integer, then $K_{1,\frac{1}{\rho}-1}$ shows that the bound in Theorem \ref{theorem2} is tight.

\end{document}